# SINGULARITY POINTS FOR FIRST PASSAGE PERCOLATION

By J. E. Yukich[1] and Yu Zhang[2]

*Lehigh University and University of Colorado*

Let $0 < a < b < \infty$ be fixed scalars. Assign independently to each edge in the lattice $\mathbb{Z}^2$ the value $a$ with probability $p$ or the value $b$ with probability $1-p$. For all $u, v \in \mathbb{Z}^2$, let $T(u,v)$ denote the first passage time between $u$ and $v$. We show that there are points $x \in \mathbb{R}^2$ such that the "time constant" in the direction of $x$, namely, $\lim_{n \to \infty} n^{-1} \mathbf{E}_p[T(\mathbf{0}, nx)]$, is not a three times differentiable function of $p$.

**1. Introduction, main results.** Consider the following simple model of first passage percolation. $E := E(\mathbb{Z}^2)$ denotes the edges in the integer lattice $\mathbb{Z}^2$, $0 < a < b < \infty$ are fixed scalars, and $\Omega := \{a, b\}^E$. For all $e \in E$ and $\omega_e \in \Omega$, $P[\omega_e = a] = p$ and $P[\omega_e = b] = 1 - p$, where $0 < p < 1$. In other words, we assign either $a$ or $b$ to each edge with probability $p$ or $1 - p$ independently from the other edges. Denote the product measure on $\Omega$ by $\mathbf{P}_p$ and the expectation with respect to $\mathbf{P}_p$ by $\mathbf{E}_p$.

For all $u, v \in \mathbb{Z}^2$, let $T(u,v)$ denote the first passage time between $u$ and $v$. Formally, $T(u,v)$ is the infimum of $\sum_{e \in \gamma} w_e$, where $\gamma$ ranges over all finite paths in $\mathbb{Z}^2$ from $u$ to $v$. If $x$ and $y$ are in $\mathbb{R}^2$, we define $T(x, y) = T(x', y')$, where $x'$ (resp. $y'$) is the point in $\mathbb{Z}^2$ closest to $x$ (resp. $y$). Any possible ambiguity can be avoided by ordering $\mathbb{Z}^2$ and taking the point in $\mathbb{Z}^2$ smallest for this order.

Let $\mathbf{0}$ denote the origin of $\mathbb{R}^2$ and for all $x \in \mathbb{R}^2$, let $T(x) := T(\mathbf{0}, x)$ be the first passage time between $\mathbf{0}$ and $x$. It is well known by Kingman's subadditive ergodic theorem ((1.13) of [9]) that, for all $x \in \mathbb{R}^2$, there is a constant $\mu_p(x)$, such that

$$(1.1) \qquad \lim_{n \to \infty} \frac{T(nx)}{n} = \mu_p(x) \qquad \text{a.s. and in } L^1.$$

Received September 2004; revised June 2005.
[1]Supported in part by NSA Grant MDA904-01-1-0029 and NSF Grant DMS-02-03720.
[2]Supported in part by NSF Grant DMS-04-05150.

When $x = (1,0)$, the limit $\mu_p^* := \mu_p((1,0))$ is called the *time constant* of Hammersley and Welsh [8]. Without loss generality, for any $x \in \mathbb{R}^2$, we also call $\mu_p(x)$ the time constant in the direction of $x$.

In general, physicists believe that most percolation constants should be real analytic as functions of $p$, excepting the singularities at the critical case. In particular, when $\omega_e$ only takes value 1 or 0, the behavior of the time constant is similar to that of the correlation length [1]. Furthermore, the analyticity of the correlation length, as expected, is proved for all $p$ except for the critical case when $d = 2$ [2]. Few rigorous results are known for the time constant. Cox and Kesten (Theorem 3 of [4]) show that $\mu_p^*$ is continuous with respect to the weak convergence of the distribution of the passage times, from which it follows that $\mu_p^*$ is continuous in $p$.

With these observations, one might believe that both the correlation length and the time constant are analytic except for the critical case when $\omega_e$ takes the values 1 or 0. Furthermore, one might also expect that the behavior of the time constant in the critical case is similar to the behavior in the case when $\omega_e$ takes the values $a$ or $b$ with $0 < a < b$. We find here that the analyticity of the latter is not always true. The main goal of this paper is to show there is a direction for which the directional asymptotic speed is not three times differentiable in the parameter $p$.

Recall that the classical grid $\mathcal{L}$ for oriented percolation is given by $\mathcal{L} := \{(m,n) \in \mathbb{Z}^2 : m+n \text{ has even parity}, n \geq 0\}$. Thus, $\mathcal{L}$ is $\mathbb{Z}^2$ rotated by $\pi/4$ and correctly dilated. Let $E(\mathcal{L})$ be the edges from $(m,n) \in \mathcal{L}$ to $(m+1, n+1)$ and to $(m-1, n+1)$. To each edge $e \in E(\mathcal{L})$, we assign a passage time $a > 0$ with probability $p$ and a time $b > a$ with probability $1-p$. Henceforth, let $\Omega := \{a,b\}^{E(\mathcal{L})}$.

Let $\overrightarrow{p_c}$ denote the critical probability for oriented Bernoulli percolation on $\mathcal{L}$. For all $p \in (\overrightarrow{p_c}, 1]$, consider all paths starting from $\{(x,y) \subset \mathbb{Z}^2 : x \leq 0, y = 0\}$ in the *oriented* graph using $n$ type $a$ oriented edges $E(\mathcal{L})$ and let $(r_n(p), n)$ denote the rightmost point ("right-hand edge") of all such paths. We will often simply refer to the scalar $r_n(p)$ as the right-hand edge. In the super-critical regime $p \in (\overrightarrow{p_c}, 1]$, the rightmost point $(r_n(p), n)$ satisfies

$$\text{(1.2)} \qquad \lim_{n \to \infty} \frac{r_n(p)}{n} = \alpha(p) \qquad \text{a.s. and in } L^1,$$

as well as a central limit theorem [10]. Here $\alpha(p) \in (0,1]$ is called the asymptotic speed of super-critical oriented percolation on the edges of $\mathcal{L}$. It describes the drift of the rightmost point at level $n$.

If $p > \overrightarrow{p_c}$, then the asymptotic shape [the unit radius ball for the norm induced by the map $x \to \mu_p(x)$] exhibits a flat edge [6], which is related directly to the possibility of percolating with edges having passage time $a$. The flat edges of the asymptotic shaper are in the coordinate directions and are described analytically by Marchand [12] (see especially Theorem 1.3).

Let $p_0 \in (\vec{p_c}, 1)$ be fixed. For all $p \in (\vec{p_c}, 1)$, define a *time constant* in the direction of the critical vector with components $\alpha(p_0)$ and 1, that is, set

$$f_{p_0}(p) := \lim_{n \to \infty} \frac{\mathbf{E}_p[T((\alpha(p_0)n, n))]}{n}.$$

It is easy to see (cf. Lemma 3.3 below for details) that if $p \geq p_0$, then on the average there is an oriented path between $\mathbf{0}$ and $(\alpha(p_0)n, n)$ consisting of edges having passage time $a$, that is, $f_{p_0}(p) = a$ for all $p \in [p_0, 1]$. Thus, if $p \longmapsto f_{p_0}(p)$ is three times differentiable at $p = p_0$, then the third derivative must be zero. However, in what follows, we show there is a constant $C > 0$ such that, for all $p \in (\vec{p_c}, p_0)$, we have

(1.3) $$f_{p_0}(p) \geq a + C(p_0 - p)^2 / (-\log(p_0 - p)).$$

This is enough to show that $p \longmapsto f_{p_0}(p)$ is not three times differentiable at $p_0$. This is our main result, formally stated as follows:

THEOREM 1.1. *For all $p_0 \in (\vec{p_c}, 1)$, the function $p \longmapsto f_{p_0}(p)$ is not three times differentiable at $p = p_0$.*

REMARKS. 1. Hammersley and Welsh conjecture (Corollary 6.5.5 of [8]) that $\mu_p^*$ is concave in $p$ and thus differentiable for almost all $p$. One might also expect that $p \longmapsto f_{p_0}(p)$ is concave and differentiable, but we are unable to show it.

2. Theorem 1.1 can be generalized to include passage times having a common distribution $p\delta_a + (1-p)U(b)$, where $0 < a < b$, $p \in [0, 1]$, and $U(b)$ is an independent random variable bounded below by $b$. It is unclear (at least to us) whether Theorem 1.1 remains true for (i) more general passage times, or (ii) directions other than $(\alpha(p_0)n, n)$. It is also unclear whether the lower bound (1.3) can be improved to $f_{p_0}(p) \geq a + C(p_0 - p)/(-\log(p_0 - p))$.

3. A natural problem involves studying the properties of the asymptotic shape at the end of its flat edge for a fixed $p$. Our methods do not yield any information here.

**2. Probability bounds for the right-hand edge of super-critical percolation.** The following proposition is of independent interest and provides exponential tail bounds for the right-hand edge $r_n(p), p \in (\vec{p_c}, 1]$. We will make critical use of this estimate in the sequel, but for now we note that Proposition 2.1 should be compared with the general tail bounds of Kuczek and Crank [11] (Theorem 1, part 1), who show, for all $p \in (\vec{p_c}, 1]$ and all $0 < \varepsilon < 1$, that there are constants $K_1 := K_1(p, \varepsilon)$ and $K_2 := K_2(p, \varepsilon)$ such that, for all $n = 1, 2, \ldots,$

$$\mathbf{P}_p[r_n(p) \geq (\alpha(p) + \varepsilon)n] \leq K_1 n^{-1/2} \exp(-K_2 n).$$

PROPOSITION 2.1. *For all $q \in (\vec{p_c}, 1]$, there exists $C_1 := C_1(q) > 0$ such that for all $0 < \varepsilon < 1$, all $p \in [q, 1]$, and all $n = 1, 2, \ldots$,*

$$\mathbf{P}_p[r_n(p) \geq (\alpha(p) + \varepsilon)n] \leq C_1 n \exp(-\varepsilon^2 n / C_1).$$

The proof of Proposition 2.1 involves consideration of the renewal process arising by breaking the behavior of the rightmost point $r_n(p)$ into independent pieces, an approach developed by Kuczek [10]. Our methods require an exponential decay result on the size of a finite cluster in super-critical oriented percolation [5].

Before proving Proposition 2.1, we require some terminology [10] and a lemma. Given vertices $u$ and $v$ in $\mathcal{L}$, we say $u \to v$ if there is a sequence $v_0 = u, v_1, \ldots, v_m = v$ of points of $\mathcal{L}$ with $v_i := (x_i, y_i)$ and $v_{i+1} := (x_{i+1}, y_i + 1)$ for $0 \leq i \leq m - 1$ such that $v_i$ and $v_{i+1}$ are connected by an edge with weight $a$. Thus, $u \to v$ if there is a sequence of oriented edges each with weight $a$ joining $u$ to $v$. For $A \subset \mathbb{Z}$, let

$$\xi_n^A := \{x : (x, n) \in \mathcal{L} \text{ and } \exists x' \in A \text{ such that } (x', 0) \to (x, n) \text{ for } n > 0\}.$$

As in [10], denote the event that there exists an infinite oriented path of $a$ edges starting from $(x, y)$ by $\Omega_\infty^{(x,y)}$. We let $\xi_0' := \xi_0^{(0,0)} := \{\mathbf{0}\}$ and set

$$\xi_1' := \begin{cases} \xi_1^{(0,0)}, & \text{if } \xi_1^{(0,0)} \neq \varnothing, \\ \{1\}, & \text{otherwise,} \end{cases}$$

and define inductively, for all $n = 1, 2, \ldots$,

$$\xi_{n+1}' := \begin{cases} \{x : (x, n+1) \in \mathcal{L} \text{ and} \\ \quad (y, n) \to (x, n+1) \text{ for some } y \in \xi_n'\}, & \text{if this set is nonempty,} \\ \{n+1\}, & \text{otherwise.} \end{cases}$$

We have suppressed the dependence of $\xi_n'$ on $p$ for notational convenience. Note that $\xi_n'$ is a subset of the integers between $-n$ and $n$. Let

$$r_n'(p) := \sup\{x : x \in \xi_n'\}.$$

On $\{\xi_n^{(0,0)} \neq \varnothing\}$, we have equivalence between $r_n'(p)$ and the right-hand edge $r_n(p)$. A vertex $(x, n) \in \mathcal{L}$ is said to be a *percolation point* if and only if the event $\Omega_\infty^{(x,n)}$ occurs. Let

$$T_1 := \inf\{n \geq 1 : (r_n', n) \text{ is a percolation point}\},$$
$$T_2 := \inf\{n \geq T_1 + 1 : (r_n', n) \text{ is a percolation point}\},$$
$$\vdots$$
$$T_m := \inf\{n \geq T_{m-1} + 1 : (r_n', n) \text{ is a percolation point}\},$$

where we make the convention that $\inf \varnothing = \infty$. Define

$$\tau_1 := T_1, \qquad \tau_2 := T_2 - T_1, \ldots, \tau_m := T_m - T_{m-1},$$

where $\tau_i := 0$ if $T_i$ and $T_{i-1}$ are infinite. (Note that $T_i$ and $T_{i-1}$ are finite with probability one.) Also define

$$X_1 := r'_{T_1}, \qquad X_2 := r'_{T_2} - r'_{T_1}, \ldots, X_m := r'_{T_m} - r'_{T_{m-1}},$$

where $X_i := 0$ if $T_i = \infty$ and $T_{i-1} = \infty$. The points $\{(r'_{T_i}, T_i)\}$ are called *break points* [10] since they break the behavior of the right-hand edge into i.i.d. pieces when the origin is a percolation point. Kuczek (Theorem on page 1324, [10]) proved that, conditional on $\Omega_\infty^{(0,0)}$, $\{(X_i, \tau_i)\}$ are i.i.d. with all moments. Moreover, for all $q \in (\vec{p_c}, 1]$, there exists a positive constant $C_2 := C_2(q)$ such that, for all $p \in [q, 1]$ and all $t \geq 1$,

$$(2.1) \qquad \mathbf{P}_p[\tau_1 \geq t] \leq \mathbf{P}_p[\xi_{t-1}^{(1,1)} \neq \varnothing, (1,1) \not\to \infty] \leq C_2 \exp(-t/C_2),$$

where the last inequality is as in [5], Section 12.

If we set

$$N_n := \sup\left\{m : \sum_{i=1}^m \tau_i \leq n\right\},$$

then $r_{N_n+1}$ is the location of the right-hand edge at the first "regeneration point" after time $n$. By considering $|r_{N_n+1} - r_{N_n}|$ and $|r_n - r_{N_n}|$, it easily follows that

$$(2.2) \qquad |r_{N_n+1} - r_n| \leq 2\tau_{N_n+1}$$

(see page 1331, [10] for details).

To prove Proposition 2.1, we make use of the following probability measure on $\Omega$:

$$\bar{\mathbf{P}}_p[\cdot] := \mathbf{P}_p[\cdot | \Omega_\infty^{(0,0)}].$$

Let $\bar{\mathbf{E}}_p$ denote the expected value with respect to $\bar{\mathbf{P}}_p$. If the event $\{r_n(p) \geq (\alpha(p) + \varepsilon)n\}$ occurs for a particular configuration $\omega \in \Omega$ of edges, then it also occurs for any configuration $\omega'$ whose $a$ edges are a superset of the $a$ edges in $\omega$. Thus, the event $\{r_n(p) \geq (\alpha(p) + \varepsilon)n\}$ is increasing. Similarly, $\Omega_\infty^{(0,0)}$ is an increasing event so that, by the FKG inequality,

$$\mathbf{P}_p[\Omega_\infty^{(0,0)}]\mathbf{P}_p[r_n(p) \geq (\alpha(p) + \varepsilon)n] \leq \mathbf{P}_p[r_n(p) \geq (\alpha(p) + \varepsilon)n, \Omega_\infty^{(0,0)}],$$

that is, to say,

$$\mathbf{P}_p[r_n(p) \geq (\alpha(p) + \varepsilon)n] \leq \bar{\mathbf{P}}_p[r_n(p) \geq (\alpha(p) + \varepsilon)n].$$

LEMMA 2.1. *Let $q \in (\vec{p_c}, 1]$. There exists $C_3 := C_3(q)$ such that for all $0 < \varepsilon < 1$, all $p \in [q, 1]$, and all $n = 1, 2, \ldots$,*

$$\bar{\mathbf{P}}_p[\tau_{N_n+1} \geq \varepsilon n] \leq C_3 n \exp(-\varepsilon n / C_3). \tag{2.3}$$

We defer the proof of Lemma 2.1 and instead show how it implies Proposition 2.1. For convenience, we put $\alpha := \alpha(p)$ and $r_n := r_n(p)$.

PROOF OF PROPOSITION 2.1. By the definition of $N_n$ and (2.2) we have, for all $0 < \varepsilon < 1$ and all $n = 1, 2, \ldots$,

$$\begin{aligned}
\mathbf{P}_p[r_n \geq (\alpha + \varepsilon)n] &\leq \bar{\mathbf{P}}_p[r_n \geq (\alpha + \varepsilon)n] \\
&\leq \bar{\mathbf{P}}_p[r_{N_n+1} + 2\tau_{N_n+1} \geq (\alpha + \varepsilon)n] \\
&\leq \bar{\mathbf{P}}_p[r_{N_n+1} \geq (\alpha + \varepsilon/2)n] + \bar{\mathbf{P}}_p[\tau_{N_n+1} \geq \varepsilon n / 4].
\end{aligned}$$

By Lemma 2.1 and since $\alpha \leq 1$, the above is bounded by

$$\leq \bar{\mathbf{P}}_p[X_1 + \cdots + X_{N_n+1} \geq \alpha(1 + \varepsilon/2)n] + C_3 n \exp(-\varepsilon n / 4 C_3). \tag{2.4}$$

Put $\kappa := \kappa(p) := \bar{\mathbf{E}}_p[\tau_1]$ and note that $\kappa \geq 1$ by definition of $\tau_1$. For $n \geq \kappa$, let $m := \lfloor \frac{n}{\kappa}(1 + \varepsilon/4) \rfloor$, where, for all $x \in \mathbb{R}$, $\lfloor x \rfloor$ denotes the greatest integer less than or equal to $x$. It follows that the above is less than or equal to

$$\sum_{i=1}^{m} \bar{\mathbf{P}}_p[X_1 + \cdots + X_i \geq \alpha(1 + \varepsilon/2)n] + \bar{\mathbf{P}}_p[N_n + 1 \geq m + 1]$$
$$+ C_3 n \exp(-\varepsilon n / 4 C_3).$$

Denote the first two terms in the above inequality by $I$ and $II$. For simplicity, we put $Y_j := \kappa - \tau_j$. Thus, by definition of $\kappa$,

$$\begin{aligned}
II := \bar{\mathbf{P}}_p[N_n + 1 \geq m + 1] &= \bar{\mathbf{P}}_p\left[\sum_{j=1}^{m} \tau_j \leq n\right] = \bar{\mathbf{P}}_p\left[\sum_{j=1}^{m}(\kappa - Y_j) \leq n\right] \\
&\leq \bar{\mathbf{P}}_p\left[\sum_{j=1}^{m} Y_j \geq \kappa(n/\kappa + \varepsilon n / 4\kappa - 1) - n\right] \\
&= \bar{\mathbf{P}}_p\left[\sum_{j=1}^{m} Y_j + \kappa \geq \varepsilon n / 4\right].
\end{aligned}$$

By Markov's inequality, for all $r > 0$,

$$II \leq \exp(r\kappa) \exp(-r\varepsilon n / 4) \bar{\mathbf{E}}_p \exp\left(r \sum_{j=1}^{m} Y_j\right). \tag{2.5}$$

Since $\bar{\mathbf{E}}_p[Y_1] = 0$ and since all moments of $Y_1$ exist, it follows that, for all $p \in [q,1]$, there exists $C_4 := C_4(q)$ such that $\log \bar{\mathbf{E}}_p[\exp(rY_1)] \le C_4 r^2$ if $r < r_0 := r_0(q)$. Thus, for $r < r_0(q)$, we obtain

$$II \le \exp(r\kappa - r\varepsilon n/4 + C_4 m r^2).$$

If we let $r := \varepsilon\kappa/C$ and increase $C$ if necessary, then it follows that there exists $C_5 := C_5(q)$ such that, for all $0 < \varepsilon < 1$, all $n \ge \kappa$ and $p \in [q,1]$,

(2.6) $$II \le C_5 \exp(-\varepsilon^2 n/C_5).$$

Increasing the value of $C_5$ if necessary, we see that (2.6) holds for $n \in [1,\kappa]$ as well.

Now we bound term $I$. By Lemma 1 of [13], we know $\alpha = \bar{\mathbf{E}}_p X_1 / \kappa$ and thus, by definition of $m$, we have, for all $1 \le i \le m$,

$$\bar{\mathbf{E}}_p[X_1 + \cdots + X_i] = i\bar{\mathbf{E}}_p X_1 \le n \frac{\bar{\mathbf{E}}_p X_1}{\kappa}(1 + \varepsilon/4)$$
$$= \alpha n (1 + \varepsilon/4).$$

Thus,

$$I \le \sum_{i=1}^m \bar{\mathbf{P}}_p\left[\sum_{j=1}^i (X_j - \bar{\mathbf{E}}_p X_j) \ge \alpha n(1 + \varepsilon/2) - \alpha n(1 + \varepsilon/4)\right]$$
$$= \sum_{i=1}^m \bar{\mathbf{P}}_p\left[\sum_{j=1}^i (X_j - \bar{\mathbf{E}}_p X_j) \ge \alpha \varepsilon n/4\right].$$

Since $|X_j| \le |\tau_j|$ for all $j \le i$, where $i \le m \le 2n$, we may follow the approach used for the bound (2.6) to conclude that there exists $C_6 := C_6(q)$ such that, for all $0 < \varepsilon < 1$, $p \in [q,1]$, and all $n = 1, 2, \ldots$,

(2.7) $$I \le C_6 n \exp(-\varepsilon^2 n/C_6).$$

Recalling that

$$\mathbf{P}_p[r_n \ge (\alpha + \varepsilon)n] \le I + II + C_3 n \exp(-\varepsilon n/4C_3)$$

and applying the bounds (2.6) and (2.7), we obtain Proposition 2.1 as desired. $\square$

Now it remains to show Lemma 2.1.

PROOF OF LEMMA 2.1. By definition of $N_n$, we have, for all $0 < \varepsilon < 1$, all $p \in (\vec{p_c}, 1]$, and all $n = 1, 2, \ldots$,

$$\bar{\mathbf{P}}_p[\tau_{N_n+1} \ge \varepsilon n] = \sum_{i=1}^\infty \bar{\mathbf{P}}_p[\tau_{i+1} \ge \varepsilon n, N_n = i]$$

$$= \sum_{i=1}^{\infty} \bar{\mathbf{P}}_p \left[ \tau_{i+1} \geq \varepsilon n, \sum_{k=1}^{i} \tau_k \leq n, \sum_{k=1}^{i+1} \tau_k > n \right]$$

$$= \sum_{j \geq \varepsilon n} \sum_{i=1}^{\infty} \bar{\mathbf{P}}_p \left[ \tau_{i+1} = j, \sum_{k=1}^{i} \tau_k \leq n, \sum_{k=1}^{i} \tau_k > n - j \right].$$

Under the measure $\bar{\mathbf{P}}_p$, the $\{\tau_i\}$ are independent and, thus, the above equals

$$\sum_{j \geq \varepsilon n} \bar{\mathbf{P}}_p[\tau_{i+1} = j] \sum_{i=1}^{\infty} \bar{\mathbf{P}}_p \left[ \sum_{k=1}^{i} \tau_k \leq n, \sum_{k=1}^{i} \tau_k > n - j \right]$$

$$\leq \sum_{j \geq \varepsilon n} \bar{\mathbf{P}}_p[\tau_{i+1} = j] \sum_{i \leq 2n/\kappa} \bar{\mathbf{P}}_p \left[ \sum_{k=1}^{i} \tau_k \leq n, \sum_{k=1}^{i} \tau_k > n - j \right]$$

$$+ \sum_{i > 2n/\kappa} \bar{\mathbf{P}}_p \left[ \sum_{k=1}^{i} \tau_k \leq n, \sum_{k=1}^{i} \tau_k > n - j \right]$$

$$:= I + II.$$

Let us bound $II$. Notice that if $i > 2n/\kappa$, then $i\kappa - n > i\kappa/2$, so we have

$$\bar{\mathbf{P}}_p \left[ \sum_{k=1}^{i} \tau_k \leq n \right] = \bar{\mathbf{P}}_p \left[ \sum_{k=1}^{i} (\kappa - \tau_k) \geq i\kappa - n \right] \leq \bar{\mathbf{P}}_p \left[ \sum_{k=1}^{i} (\kappa - \tau_k) \geq \frac{i\kappa}{2} \right].$$

By the methods used to obtain (2.6), there exists $C_7 := C_7(q)$ and $C_8 := C_8(q)$ such that, for all $p \in [q, 1]$ and all $n = 1, 2, \ldots$,

(2.8) $$II \leq \sum_{i \geq n/\kappa + n} C_7 \exp(-i/C_7) \leq C_8 \exp(-n/C_8).$$

Let us bound term $I$. The second factor in $I$ is bounded by the number of summands showing that

$$I \leq \left( \frac{2n}{\kappa} \right) \sum_{j \geq \varepsilon n} \bar{\mathbf{P}}_p[\tau_1 = j] \leq 2n \bar{\mathbf{P}}_p[\tau_1 \geq \varepsilon n],$$

since $\kappa \geq 1$. Combining this with (2.1) shows that there exists $C_9 := C_9(q)$ such that, for all $0 < \varepsilon < 1$, all $p \in [q, 1]$, and all $n = 1, 2, \ldots$,

$$I \leq C_9 n \exp(-\varepsilon n/C_9).$$

Lemma 2.1 now follows from (2.8) and the above inequality. □

**3. Auxiliary lemmas.** The proof of Theorem 1.1 rests on the upper bound for the right-hand edge of supercritical percolation (Proposition 2.1), as well as a lower bound for first passage times, given in the upcoming Proposition 4.1. Before proving the latter, we require six straightforward lemmas. Our first lemma gives a way to prove the asserted nondifferentiability of $f_{p_0}$, where we recall that $p_0 \in (\vec{p_c}, 1)$ is fixed once and for all. Let log denote the natural logarithm. For the remainder of the paper, we fix $q \in (\vec{p_c}, p_0)$.

LEMMA 3.1. *Suppose $h:[0,1] \to \mathbb{R}^+$ satisfies $h(p) = 0$ for all $p \geq p_0$. If there exists $\delta := \delta(q) > 0$ such that, for all $p \in [q, p_0)$,*

$$h(p) \geq \frac{\delta(p_0 - p)^2}{\log(1/(p_0 - p))}, \tag{3.1}$$

*then $h'''(p_0)$ does not exist.*

PROOF. We use elementary calculus. If $h'''(p_0)$ did exist, then necessarily $h'''(p_0) = h''(p_0) = h'(p_0) = 0$. It follows that $|h''(p)| = |h''(p) - h''(p_0)| \leq |p_0 - p|$ if $|p - p_0|$ is small enough. For such $p$, we have $|h'(p)| = |\int_{p_0}^p h''(u)\,du| \leq \int_p^{p_0} |h''(u)|\,du \leq (p_0 - p)^2$, that is, $h'(p)$ grows at most like a quadratic in $p_0 - p$. Similarly, $h(p)$ grows at most like a cubic in $p_0 - p$ for $|p - p_0|$ small enough. This is a contradiction. □

To show that the function $f_{p_0}$ of Theorem 1.1 satisfies the conditions of Lemma 3.1, we will need several more lemmas and a proposition.

LEMMA 3.2. *For all $p \in (\vec{p_c}, p_0]$, we have $\alpha(p_0) - \alpha(p) \geq 2(p_0 - p)$.*

PROOF. See [5], page 1006, display (12). □

LEMMA 3.3. $f_{p_0}(p) = a$ *for all $p \in [p_0, 1]$.*

PROOF. By the central limit theorem of Kuczek (Corollary 1 of [10]), with probability $1 - o(1)$, there exists an oriented path $\gamma$ of $n$ type $a$ edges, starting at **0** and terminating at a point $(x, n)$, where $\alpha(p_0)n < x$. Similarly, reversing the orientation of the edges, with probability $1 - o(1)$, there exists a path $\gamma'$ of $n$ type $a$ oriented edges, starting at $(\alpha(p_0)n, n)$ and terminating at a point $(s, 0)$, where $s \geq \alpha(p_0)n$. The paths $\gamma$ and $\gamma'$ intersect at some point $Q \in \mathbb{Z}^2$. Let $\gamma_1$ be the restriction of $\gamma$ between **0** and $Q$; let $\gamma_1'$ be the restriction of $\gamma'$ between $Q$ and $(\alpha(p_0)n, n)$. Let $\gamma_u$ be the union of $\gamma_1$ and $\gamma_1'$. Then $\gamma_u$ is an oriented path $\mathbf{0} \to Q \to (\alpha(p_0)n, n)$ consisting exclusively of $n$ type $a$ edges showing that

$$T((\alpha(p_0)n, n)) = an \tag{3.2}$$

on a set with probability $1 - o(1)$. Since $n^{-1}T((\alpha(p_0)n, n))$ is bounded by $b$, the conclusion follows. $\square$

We will adhere to the following terminology throughout. Given a path $\gamma$ in the lattice $\mathcal{L}$, $T(\gamma)$ denotes its weight $\sum_{e \in \gamma} \omega_e$, where $P[\omega_e = a] = p, P[\omega_e = b] = 1 - p$. We let $\mathcal{P}(\alpha(p_0)n)$ denote all paths (oriented or not) $\gamma: \mathbf{0} \mapsto (\alpha(p_0)n, n)$ *in the lattice* $\mathcal{L}$ whose weight equals the first passage time $T((\alpha(p_0)n, n))$. [If $x \in \mathbb{R}$, then we adopt the convention that the path $\gamma: \mathbf{0} \longmapsto (x, n)$ denotes the path between $\mathbf{0}$ and $(\lfloor x \rfloor, n)$.] If $p \in (\vec{p_c}, p_0]$, then $T(\gamma), \gamma \in \mathcal{P}(\alpha(p_0)n)$, will tend to exceed $an$, since typically, under $\mathbf{P}_p$, the edges in $\gamma$ required to link $\mathbf{0}$ with points to the right of $(\alpha(p)n, n)$, for example, $(\alpha(p_0)n, n)$, will not all have weight $a$.

Consider $\delta := \delta(q) \in (0, 1/2)$ with a value to be specified later. For all $p \in [q, p_0]$, let $\mathcal{P}_n := \mathcal{P}_n(p_0, p, \delta) \subset \mathcal{P}(\alpha(p_0)n)$ be the (possibly empty) subset of $\mathcal{P}(\alpha(p_0)n)$ consisting of paths $\gamma$ whose weight satisfies

$$T(\gamma) \leq an\left(1 + \frac{\delta(p_0 - p)^2}{\log(1/(p_0 - p))}\right).$$

Thus, $\mathcal{P}_n \neq \varnothing$ is the event that the first passage time $T((\alpha(p_0)n, n))$ is bounded above by $an(1 + \frac{\delta(p_0-p)^2}{\log(1/(p_0-p))})$. We will show in Proposition 4.1 below that the probability of $\mathcal{P}_n \neq \varnothing$ is exponentially small, but first we require a few more lemmas. Recalling that $\vec{p_c} < q < p_0 < 1$ and $p \in [q, p_0]$, we will henceforth assume, without loss of generality, that $q$ is close enough to $p_0$ to guarantee that

(3.3) $$\frac{a\nu}{\log(1/(p_0 - p))} \leq 1 \quad \text{and} \quad \log\left(\frac{1}{p_0 - p}\right) > 1.$$

LEMMA 3.4. *If $\gamma \in \mathcal{P}_n$, then $\gamma \subset [-2n, 2n] \times [-n, 2n]$.*

PROOF. It suffices to show that if $\gamma \in \mathcal{P}_n$, then $\gamma$ has at most $2n$ edges. Since $\delta < 1/2$ and $\frac{(p_0-p)^2}{\log(1/(p_0-p))} < 1$, it follows that if $\gamma \in \mathcal{P}_n$, then $T(\gamma) < 2an$. Since every edge in $\gamma$ has weight at least $a$, it follows that $\gamma$ has at most $2n$ edges. $\square$

Given $\gamma \in \mathcal{P}(\alpha(p_0)n)$, an edge $e := ((x_1, y_1), (x_2, y_2))$ belonging to $\gamma$ is termed "repeated" if the horizontal strip $\mathbb{R} \times [y_1, y_2]$ contains at least one other edge in $\gamma$ and to the left of $e$. Edges $e \in \gamma$ are called "sub-optimal" if either $e$ has weight $b$ or if $e$ is repeated. Roughly speaking, paths $\gamma \in \mathcal{P}_n$ cannot use many sub-optimal edges. Edges $e := (u, v)$ are considered to be closed line segments in $\mathbb{R}^2$ in the sense that $e$ contains its endpoints $\{u\}$ and $\{v\}$.

LEMMA 3.5. *Let $\nu := (\min(b-a, a))^{-1}$. If $\gamma \in \mathcal{P}_n$, then there are at most*

$$(3.4) \qquad k := k(p, p_0, n) := \left\lfloor \frac{a\nu\delta(p_0 - p)^2 n}{\log(1/(p_0 - p))} \right\rfloor$$

*sub-optimal edges in $\gamma$.*

PROOF. Each sub-optimal edge in $\gamma$ contributes an extra cost of at least $\min(b-a, a)$. □

Recalling that $\vec{p_c} < q < p_0 < 1$ and $p \in [q, p_0]$, we will henceforth assume, without loss of generality, that $q$ is close enough to $p_0$ to guarantee that (3.3) holds and that $k \in [0, \frac{n}{10}]$. Given $\gamma \in \mathcal{P}_n$, project all sub-optimal edges in $\gamma$ onto the $x$-axis. The projection forms a possibly empty collection of closed intervals on the $x$-axis which may overlap. However, when the projection is nonempty, the union forms a collection of *disjoint* closed intervals $I_1(\gamma), I_2(\gamma), \ldots, I_j(\gamma)$ called the $x$-trace $\tau_x(\gamma)$ of $\gamma \in \mathcal{P}_n$. The intervals in $\tau_x(\gamma)$ have integral endpoints and belong to $[-2n, 2n]$ by Lemma 3.4. Here $j \in \mathbb{N}$ cannot exceed the number $k$ of sub-optimal edges; if $k = 0$, then there is no $x$-trace. Note that distinct paths $\gamma \in \mathcal{P}_n$ may have identical $x$-traces.

DEFINITION 3.1. For all $1 \leq j \leq k$, let $\mathcal{T}_j^x$ denote the collection of all $x$-traces consisting of $j$ disjoint subintervals.

Next, given $\gamma \in \mathcal{P}_n$, remove all edges in $\gamma$ whose projection onto the $x$-axis is a proper subset of $\tau_x(\gamma)$ (some such edges may be oriented and have weight $a$). What remains are called the *optimal* edges in $\gamma$; such edges are necessarily oriented up edges with weight $a$. By definition, these edges collectively form a sequence of disjoint paths $\gamma_1, \gamma_2, \ldots$, each consisting of *oriented* edges having weight $a$. We call $\gamma_1, \gamma_2, \ldots$, "*optimal paths.*" Note that optimal paths lie in $[-2n, 2n] \times [0, n]$.

Observe that the $\gamma_i, i \geq 1$, are contained in the horizontal strips $\mathbb{R} \times [y_i, y_i']$, where $y_i$ and $y_i'$ denote the $y$ coordinates of the initial and terminal points of $\gamma_i$, respectively.

We project all optimal edges in $\gamma$ onto the (vertical) $y$-axis. The projection yields a collection of intervals $I_1'(\gamma), I_2'(\gamma), \ldots$, which we call the $y$-trace $\tau_y(\gamma)$ of $\gamma$. Each interval in $\tau_y(\gamma)$ is a subset of $[0, n]$.

DEFINITION 3.2. For all $1 \leq j \leq k$, let $\mathcal{T}_j^y$ denote the collection of all $y$-traces consisting of $j$ subintervals.

Given $\gamma \in \mathcal{P}_n$, we call the set of intervals $\tau_{xy} := \{I_i(\gamma)\}_{i=1}^{j_1} \cup \{I_i'(\gamma)\}_{i=1}^{j_2}$ the $xy$-trace of $\gamma$. The collection of $xy$-traces will provide a convenient combinatorial way to upper bound the probability that $\mathcal{P}_n \neq \varnothing$. Since the number of optimal paths differs from the number of disjoint intervals in the $x$-trace by at most one, it follows that $|j_1 - j_2| \leq 1$. We say that $\tau_{xy}$ is an $xy$-trace of cardinality $j$ if $j_1 \vee j_2 = j$. Considering the three cases $j_1 = j_2, j_1 = j_2 - 1$, and $j_2 = j_1 - 1$, we see that the collection of all $xy$-traces of cardinality $j$ has the representation

$$\mathcal{T}_j := \{(I_i, I'_i)_{i=1}^j : I_i \in \mathcal{T}_j^x, I'_i \in \mathcal{T}_j^y\}$$

$$\cup \{(I_i, I'_i)_{i=1}^j : I_i \in \mathcal{T}_{j-1}^x, I_j = \varnothing, I'_i \in \mathcal{T}_j^y\}$$

$$\cup \{(I_i, I'_i)_{i=1}^j : I_i \in \mathcal{T}_j^x, I'_i \in \mathcal{T}_{j-1}^y, I'_j = \varnothing\}.$$

Since elements of $\mathcal{T}_j^x$ and $\mathcal{T}_j^y$ have integral endpoints, Lemma 3.4 implies that card $\mathcal{T}_j^x \leq \binom{4n}{2j}$. Notice that the elements of $\mathcal{T}_j^y$ have integral endpoints which may coincide (they coincide if there is an integer $i$ such that $y'_i = y_{i+1}$). The elements of $\mathcal{T}_j^y$ can be coded by their endpoints $\{(y_i, y'_i)\}_{i=1}^j$, so that, for example, the sequence $1, 2, 2, 5, 7, 8$ denotes the following three intervals on the $y$-axis: $I'_1 := ((0,1), (0,2)), I'_2 := ((0,2), (0,5)), I'_3 := ((0,7), (0,8))$. Clearly, $\mathcal{T}_j^y \leq \binom{2n}{2j}$. Since clearly $\binom{2n}{2j} \leq \binom{4n}{2j}$ for $1 \leq j \leq k$, we deduce the crude bound:

LEMMA 3.6. *For all $1 \leq j \leq k$, we have* card $\mathcal{T}_j \leq 3\binom{4n}{2j}^2$.

**4. Lower bounds for first passage times.** Recall that $q$ and $p_0$ are fixed scalars satisfying $\vec{p_c} < q < p_0$. By Lemma 3.3, we have $f_{p_0}(p) - a = 0$ for all $p \in [p_0, 1]$. It remains to show that $f_{p_0} - a$ satisfies inequality (3.1). We do this by showing that the first passage time $T((\alpha(p_0)n, n))$ is bounded below by

$$an\left(1 + \frac{\delta(p_0 - p)^2}{\log(1/(p_0 - p))}\right),$$

with overwhelming probability for $p \in [q, p_0]$. Recalling the definition of $C_1$ in Proposition 2.1, we have the following:

PROPOSITION 4.1. *For all $p \in [q, p_0]$ and all $n = 1, 2, \ldots$,*

$$\mathbf{P}_p[\mathcal{P}_n(p_0, p, \delta) \neq \varnothing] \leq C_1 n^2 \exp(-(p_0 - p)^2 n/4C_1).$$

Before proving Proposition 4.1, we first show how it implies that $f_{p_0} - a$ satisfies the conditions of Lemma 3.1. We have, for all $p \in [q, p_0]$,

$$f_{p_0}(p) := \lim_{n \to \infty} \frac{\mathbf{E}_p[T((\alpha(p_0)n, n))]}{n}$$

$$\geq \liminf_{n \to \infty} \frac{\mathbf{E}_p[T((\alpha(p_0)n, n))\mathbb{1}_{\mathcal{P}_n = \varnothing}]}{n}$$

$$\geq a + \frac{a\delta(p_0 - p)^2}{\log(1/(p_0 - p))}$$

by Proposition 4.1 and since $T((\alpha(p_0)n, n)) \leq bn$. Since $\delta > 0$, then together with Lemma 3.3, this shows that $f_{p_0} - a$ satisfies the conditions of Lemma 3.1, concluding the proof of Theorem 1.1.

Roughly speaking, Proposition 4.1 holds for the following reasons. If $T((\alpha(p_0)n, n))$ is small [i.e., bounded above by $an(1 + \frac{\delta(p_0-p)^2}{\log(1/(p_0-p))})$], then the shortest travel time path cannot have too many sub-optimal edges. The path to $(\alpha(p_0)n, n)$ is thus nearly an oriented path with only $a$ edges. However, with such edges, an oriented path will typically only reach $(\alpha(p)n, n)$, where $\alpha(p) < \alpha(p_0)$. The estimate of the probability of the complement of such an event is handled by Proposition 2.1 and some combinatorial estimates.

We note here that if $T((\alpha(p_0)n, n))$ could be bounded above by $an(1 + \frac{\delta(p_0-p)}{\log(1/(p_0-p))})$ with high probability, then our proof would show that $p \mapsto f_{p_0}(p)$ is not two times differentiable at $p = p_0$. We are unfortunately unable to show such a bound.

To prove Proposition 4.1, we introduce some terminology. Given $l = 1, 2, \ldots$, say that a path $\gamma$ has rightward displacement of $l$ if the difference between the $x$-components of the terminal and initial points of $\gamma$ equals $l$. For all integral $m \in [n-k, n]$, $\varepsilon > 0$, and $p \in [q, 1]$, let $D(n, m, p, \varepsilon) \subset \Omega$ denote the event that there exists an optimal path beginning at $\mathbf{0}$ containing $m$ edges, and with rightward displacement at least $(\alpha(p) + \varepsilon)n$. Proposition 2.1 implies, for all $p \in [q, 1]$ and all $n = 1, 2, \ldots$,

$$\begin{aligned}
\mathbf{P}_p[D(n, m, p, \varepsilon)] &\leq \mathbf{P}_p[r_m \geq (\alpha(p) + \varepsilon)n] \\
&\leq C_1 m \exp(-\varepsilon^2 m / C_1) \\
&\leq C_1 n \exp(-\varepsilon^2 n / 2C_1)
\end{aligned} \quad (4.1)$$

since $\frac{9n}{10} \leq m \leq n$. We are now ready to provide the following:

PROOF OF PROPOSITION 4.1. Let $p \in [q, p_0]$ and suppose $\mathcal{P}_n \neq \varnothing$. For any $\gamma \in \mathcal{P}_n$, let $d_{\mathrm{opt}}(\gamma)$ be the total rightward displacement by the optimal edges in $\gamma$. In other words, $d_{\mathrm{opt}}(\gamma)$ is the combined length of the projection of the optimal edges in $\gamma$ onto the $x$-axis. Equivalently, $d_{\mathrm{opt}}(\gamma)$ is the difference between the rightward displacement of $\gamma$ and the sum of the lengths of the intervals in the $x$-trace $\tau_x(\gamma)$. For any $\gamma \in \mathcal{P}_n$, we clearly have $d_{\mathrm{opt}}(\gamma) \geq \alpha(p_0)n - k$, that is,

$$d_{\mathrm{opt}}(\gamma) \geq \alpha(p_0)n - \left\lfloor \frac{a\nu\delta(p_0 - p)^2 n}{\log(1/(p_0 - p))} \right\rfloor$$

$$\geq \alpha(p)n + \left(\frac{\alpha(p_0) - \alpha(p)}{2}\right)n$$
$$+ \left\{\left(\frac{\alpha(p_0) - \alpha(p)}{2}\right)n - \frac{a\nu\delta(p_0 - p)^2 n}{\log(1/(p_0 - p))}\right\}.$$

By Lemma 3.2, the term inside the braces exceeds $n(p_0 - p)(1 - \frac{a\nu\delta(p_0-p)}{\log(1/(p_0-p))})$, which by (3.3) is nonnegative. Therefore, for all $\gamma \in \mathcal{P}_n$,

$$d_{\text{opt}}(\gamma) \geq \alpha(p)n + \left(\frac{\alpha(p_0) - \alpha(p)}{2}\right)n \geq \alpha(p)n + (p_0 - p)n.$$

Let $\mathcal{P}'_n$ denote all (not necessarily oriented) paths in the lattice $\mathcal{L}$ beginning at $\mathbf{0}$ and ending at a point $(m, n), m \in \mathbb{N}$, with an $xy$-trace having cardinality at most $k$. We thus have

$$\mathbf{P}_p[\mathcal{P}_n \neq \varnothing] \leq \mathbf{P}_p[\exists \gamma \in \mathcal{P}'_n : d_{\text{opt}}(\gamma) \geq \alpha(p)n + (p_0 - p)n]$$
$$= \mathbf{P}_p[\exists \gamma \in \mathcal{P}'_n : d_{\text{opt}}(\gamma) \geq \alpha(p)n + (p_0 - p)n, \tau_{xy}(\gamma) = \varnothing]$$
$$+ \sum_{j=1}^{k} \mathbf{P}_p[\exists \gamma \in \mathcal{P}'_n : d_{\text{opt}}(\gamma) \geq \alpha(p)n + (p_0 - p)n, \tau_{xy}(\gamma) \in \mathcal{T}_j],$$

since $\mathcal{P}'_n$ is the disjoint union (over $T$ in $\mathcal{T}_j$ and $j \in \{1, 2, \ldots, k\}$) of paths in $\mathcal{L}$ beginning at $\mathbf{0}$ and having an $xy$-trace $T$ for some $T \in \mathcal{T}_j$ and some $1 \leq j \leq k$. By additivity, the above equals

$$\mathbf{P}_p[\exists \gamma \in \mathcal{P}'_n : d_{\text{opt}}(\gamma) \geq \alpha(p)n + (p_0 - p)n, \tau_{xy}(\gamma) = \varnothing]$$
(4.2)
$$+ \sum_{j=1}^{k} \sum_{T \in \mathcal{T}_j} \mathbf{P}_p[\exists \gamma \in \mathcal{P}'_n : d_{\text{opt}}(\gamma) \geq \alpha(p)n + (p_0 - p)n, \tau_{xy}(\gamma) = T].$$

Consider a fixed $xy$-trace $T \in \mathcal{T}_j$. Every such trace $T$ is uniquely defined by a set of deterministic points $\{(P_i, P'_i)\}_{i=1}^{2j}$, where $(P_i, P'_i) \in \mathcal{L}, 1 \leq i \leq 2j$, are the endpoints of $j$ optimal paths.

By independence and invariance by translation, the probability that there exists an optimal path between $(P_1, P'_1)$ and $(P_2, P'_2)$ and a second optimal path between $(P_3, P'_3)$ and $(P_4, P'_4)$ equals the probability that there exists an optimal path joining $\mathbf{0}$, the point $(P_2 - P_1, P'_2 - P'_1)$ and the point

$$((P_2 - P_1) + (P_4 - P_3), (P'_2 - P'_1) + (P'_4 - P'_3)).$$

More generally, the probability that there exist optimal paths joining $(P_i, P'_i)$ and $(P_{i+1}, P'_{i+1})$, for all $1 \leq i \leq 2j - 1$, is bounded by the probability that there exists an optimal path between $\mathbf{0}$ and $(\sum_{i=1}^{2j-1}(P_{i+1} - P_i), \sum_{i=1}^{2j-1}(P'_{i+1} - P'_i))$. Any such path has a total of $N := \sum_{i=1}^{2j-1}(P'_{i+1} - P'_i)$ edges, where $N \in [n - k, n - 1]$. Thus, for each $1 \leq j \leq k$, and each $T \in \mathcal{T}_j$, each summand in (4.2) is bounded by the probability that there is an optimal path with $N$ edges with rightward displacement at least $\alpha(p)n + (p_0 - p)n$, that is, by the probability of $D(n, N, p, p_0 - p)$. Similarly, the first probability in (4.2) is bounded by the probability of $D(n, n, p, p_0 - p)$. It follows by Lemma 3.6 and (4.1) that (4.2) becomes

$$
(4.3) \quad \mathbf{P}_p[\mathcal{P}_n \neq \varnothing] \leq C_1 n \exp\left(-\frac{(p_0 - p)^2 n}{2C_1}\right)
$$
$$
+ 3C_1 n \sum_{j=1}^{k} \binom{4n}{2j}^2 \exp\left(-\frac{(p_0 - p)^2 n}{2C_1}\right).
$$

To conclude the proof of Proposition 4.1, it suffices to show that, for all $1 \leq j \leq k$,

$$
(4.4) \quad \binom{4n}{2j}^2 \leq \exp\left(\frac{(p_0 - p)^2 n}{4C_1}\right).
$$

To do this, we will make use of ([7], Corollary 2.6.2)

$$
\binom{u}{v} \leq \exp\left(uH\left(\frac{v}{u}\right)\right), \qquad u, v \in \mathbb{N},
$$

where, for all $x \in (0, 1)$,

$$
H(x) := -x \log x - (1 - x) \log(1 - x).
$$

Thus, for all $j = 1, 2, \ldots, k := \lfloor a\nu\delta(p_0 - p)^2 n / \log(\frac{1}{p_0 - p}) \rfloor$, we have

$$
(4.5) \quad \binom{4n}{2j} \leq \binom{4n}{2k} \leq \exp\left(4nH\left(\frac{k}{2n}\right)\right),
$$

where the first inequality holds since $k \leq n/10$.

There is $x_0 \in (0, 1)$ such that if $x \in (0, x_0)$, then $-(1 - x)\log(1 - x) \leq -\log(1 - x) \leq -x \log x$, showing that, for all $x \in (0, x_0)$, we have

$$
H(x) \leq 2x \log \frac{1}{x}.
$$

By choosing $\delta := \delta(q)$ so small that $a\nu\delta < x_0$, we guarantee that $k/2n < x_0$. Since $x \log \frac{1}{x}$ is increasing on $(0, 1)$, we obtain

$$
H\left(\frac{k}{2n}\right) \leq \frac{a\nu\delta(p_0 - p)^2}{\log(1/(p_0 - p))} \log\left(\frac{2n}{\lfloor a\nu\delta(p_0 - p)^2 n / \log(1/(p_0 - p)) \rfloor}\right)
$$
$$
\leq \frac{a\nu\delta(p_0 - p)^2}{\log(1/(p_0 - p))} \log\left(\frac{4 \log(1/(p_0 - p))}{a\nu\delta(p_0 - p)^2}\right),
$$

since $\frac{x}{\lfloor y \rfloor} \leq \frac{2x}{y}$ for $x, y \geq 1$. Simple algebra shows that the above equals

$$\frac{a\nu\delta(p_0 - p)^2}{\log(1/(p_0 - p))} \left[ \log\log\left(\frac{1}{p_0 - p}\right) + \log\left(\frac{4}{a\nu\delta}\right) + 2\log\left(\frac{1}{p_0 - p}\right) \right]$$

$$< 3a\nu\delta(p_0 - p)^2 + a\nu\delta(p_0 - p)^2 \log\left(\frac{4}{a\nu\delta}\right)$$

using $-\infty < \log\log t \leq \log t$ for $t > 1$ and $\log(\frac{1}{p_0 - p}) > 1$. Choosing $\delta := \delta(q) \in (0, 1/2)$ so small that $a\nu\delta \log(\frac{4}{a\nu\delta}) \leq (a\nu\delta)^{1/2}$, we get

$$(4.6) \qquad H\left(\frac{k}{2n}\right) \leq 4(a\nu\delta)^{1/2}(p_0 - p)^2.$$

Substituting (4.6) into (4.5) and squaring, we obtain, for all $1 \leq j \leq k$,

$$\binom{4n}{2j}^2 \leq \exp(32(a\nu\delta)^{1/2}(p_0 - p)^2 n).$$

Recalling that $C_1$ depends only on $q$, we may choose $\delta := \delta(q) > 0$ even smaller if necessary to ensure that $32(a\nu\delta)^{1/2} < 1/4C_1$, thus, showing (4.4). Proposition 4.1 follows. $\square$

**Acknowledgments.** The authors gratefully acknowledge detailed referee comments, which resulted in an improved exposition. They also acknowledge the hospitality of Lehigh University and the University of Colorado, where parts of this research were done.

DEPARTMENT OF MATHEMATICS
LEHIGH UNIVERSITY
BETHLEHEM, PENNSYLVANIA 18015
USA
E-MAIL: joseph.yukich@lehigh.edu
URL: www.lehigh.edu/˜jey0/jey0.html

DEPARTMENT OF MATHEMATICS
UNIVERSITY OF COLORADO
COLORADO SPRINGS, COLORADO 80933
USA
E-MAIL: yzhang@math.uccs.edu